\documentclass[12pt]{article}

\usepackage{amsmath}
\usepackage{graphicx}
\usepackage{booktabs}
\usepackage{amsfonts}
\usepackage{hyperref}
\usepackage{amssymb}

\renewcommand{\theequation}
{\arabic{section}.\arabic{equation}}

\makeatletter
\def\eqnarray{ \stepcounter{equation} \let\@currentlabel=\theequation
 \global\@eqnswtrue
 \global\@eqcnt\z@
 \tabskip\@centering
 \let\\=\@eqncr
 $$\halign to \displaywidth\bgroup\@eqnsel\hskip\@centering
 $\displaystyle\tabskip\z@{##}$&\global\@eqcnt\@ne
 \hfil$\displaystyle{{}##{}}$\hfil
 &\global\@eqcnt\tw@$\displaystyle\tabskip\z@{##}$\hfil
 \tabskip\@centering&\llap{##}\tabskip\z@\cr}
\makeatother

\makeatletter
\def\@arrayacol{\edef\@preamble{\@preamble \hskip .5\arraycolsep}}
\def\array{\let\@acol\@arrayacol \let\@classz\@arrayclassz
\let\@classiv\@arrayclassiv \let\\\@arraycr\def\@halignto{}\@tabarray}
\makeatother

\makeatletter
\newcounter{subeqncnt}
\def\thesubeqncnt{\alph{subeqncnt}}
\def\subequations{\begingroup%
   \stepcounter{equation}\edef\@tempa{\theequation}%
   \let\c@equation\c@subeqncnt\c@subeqncnt\z@
   \edef\theequation{\@tempa\noexpand\thesubeqncnt}}

\makeatother

\newcommand{\be}{\begin{equation}}
\newcommand{\ee}{\end{equation}}

\newcommand{\bea}{\begin{eqnarray}}
\newcommand{\eea}{\end{eqnarray}}
\newcommand{\nn}{\nonumber}

\def\CR {{\cal R}}

\def\CZ {{\cal Z}}

\def\Det{{\rm Det}}

\begin{document}

\setlength{\baselineskip}{7mm}
\begin{titlepage}
 \begin{flushright}
{\tt NRCPS-HE-26-2020}
\end{flushright}

\begin{center}
{\Large ~\\{\it   Extended Kolmogorov Entropy

\vspace{1cm}

}

}

\vspace{1cm}

{\sl George Savvidy

\bigskip
\centerline{${}$ \sl Institute of Nuclear and Particle Physics}
\centerline{${}$ \sl Demokritos National Research Center, Ag. Paraskevi,  Athens, Greece}
\bigskip

}
\end{center}
\vspace{30pt}

\centerline{{\bf Abstract}}

The Kolmogorov entropy allows to split the dynamical systems that have equivalent continuous spectrum into non-isomorphic subclasses. In this paper we make an attempt to generalise the concept of entropy that will allow to split the systems with equivalent continuous spectrum and equal Kolmogorov entropies into finer non-isomorphic subclasses.  We will define and calculate the new metrical invariant for the hyperbolic automorphisms on a torus.  The hyperbolic systems on a torus are perfect pseudorandom number generators for the Monte-Carlo (MC) simulations in high energy physics, and the new metrical invariant allows a finer characterisation of the MC generators when they have equal entropies.

\vspace{12pt}

\noindent

\end{titlepage}




\pagestyle{plain}

\section{\it Introduction}
 
Hyperbolic dynamical systems $T$  defined on a torus have strong instability of their trajectories, as strong as it can be in principle \cite{anosov}. In series of publications we proposed to use the Anosov hyperbolic C-systems on a torus to generate pseudorandom numbers for the Monte-Carlo simulation  \cite{yer1986a,konstantin}.  The Kolmogorov entropy $h(T)$ \cite{kolmo,kolmo1,Shannon,sinai3} allows to quantitatively characterise the pseudorandom number generators \cite{yer1986a,konstantin} in a sense that the larger the entropies are the better the generators are \cite{konstantin,Savvidy:2015jva,Savvidy:2015ida}. At the same time we found that the random number generators with equal entropies may have essentially different stochastic properties and that there is a need for additional metrical invariant which will allow a finer characterisation of the MC generators when they have equal entropies.

As it is well known, the Kolmogorov entropy allows to classify and differentiate the dynamical systems that have equivalent continuous spectrum \cite{kolmo} (in mathematical literature: the systems with "countable Lebesgue  spectrum"). In this paper we make an attempt to generalise the concept of entropy  that will allow to differentiate the systems which have  equivalent countable Lebesgue  spectrum  \cite{Koopman} and equal Kolmogorov entropies.  That is in the cases when  the unitary operators $U_1$ and  $U_2$ associated  with the hyperbolic systems $T_1$ and $T_2$  have countable Lebesgue  spectrum and  equal entropies $h(T_1)=h(T_2)$.  It seems natural to think that the available source of information that will allow to distinguish these systems can lie in refined properties of the eigenvalues spectrum of the operator $T$.   In the next section we will define the new invariant.

\section{\it Extended Entropy}
 
The hyperbolic dynamical systems are systems that have a  uniform and exponential instability of their phase trajectories $u^{(n)}=T^n u^{(0)}$, where $T$ is an evolution operator. The exponential instability of the dynamical system  $T $ takes place when the deviation of trajectories $\delta u^{(n)} $ has an exponential character. In such systems a distance between infinitesimally close trajectories increases exponentially and on a compact phase space $u \in M$ equipped with a positive Liouville's measure $d \mu(u)$ leads to the uniform distribution of almost all trajectories over the whole phase space $M$.

For that reason the dynamical systems that have local and homogeneous hyperbolic instability of the phase trajectories  have very extended and rich ergodic properties \cite{anosov}.  As such they have 
mixing of all orders, countable Lebesgue  spectrum, positive entropy 
and occupy a nonzero volume  in the space of dynamical systems \cite{anosov}.  
The important examples of the hyperbolic C-systems\footnote{
D.V. Anosov gave the definition of C-systems in his outstanding work \cite{anosov}.
The definition of C-systems uses such mathematical concepts as a tangent vector bundle, a derivative mapping, the contracting and expanding linear spaces, foliations  and other concepts.  The review of the C-systems \cite{anosov},  of the Kolmogorov entropy \cite{kolmo,kolmo1,sinai3}, the properties of its periodic trajectories can be found in the recent article \cite{Savvidy:2015ida} and in Appendix B. } are: i)   
the geodesic flow on the Riemannian manifolds of variable negative curvature\footnote{The exponential instability of geodesics on Riemannian manifolds of constant negative curvature has been studied by Lobachevsky and Hadamard and by  Hedlund and Hopf \cite{hedlund,hopf}.} and 
ii) the hyperbolic automorphisms of tori.

Particular systems which we shall consider here are the automorphisms of a torus or of a unit hypercube  in Euclidean space $\mathbb{E}^N$ with coordinates $ (u_1,...,u_N)$ $mod 1$~\cite{anosov,yer1986a,konstantin,Savvidy:2015ida}:
\be
\label{eq:rec}
u_i^{(k+1)} = \sum_{j=1}^N T_{ij} \, u_j^{(k)} ~~~~~\textrm{mod}~ 1,~~~~~~~~~k=0,1,2,...
\ee
where the components of the vector $u^{(k)}$ are defined as 
$$
u^{(k)}= (u^{(k)}_1,...,u^{(k)}_N).
$$ 
 The phase space of the system (\ref{eq:rec} ) can also be considered as a $N$-dimensional torus appearing at  factorisation of the Euclidean space $\mathbb{E}^N$ with coordinates $u= (u_1,...,u_N)$ over an integer lattice $\CZ^N$. The dynamical system defined here by the integer matrix $T$
should have a determinant equal to one $\Det T =1$.
In order for the automorphisms (\ref{eq:rec})  to fulfil  the Anosov hyperbolicity
C-condition (see Appendix B) it is necessary
and sufficient that the matrix $T$ has no eigenvalues on the unit circle \cite{anosov}.
Therefore  the  spectrum $\{ \Lambda = {\lambda_1},...,
\lambda_N \}$ of the matrix $T$ should fulfil the following
two conditions:
\bea\label{mmatrix1}
1)~\Det  T=  {\lambda_1}\,{\lambda_2}...{\lambda_N}=1,~~~~~
2)~~\vert {\lambda_i} \vert \neq 1, ~~~\forall ~~i .~~~~~~
\eea
Because the determinant of the matrix $T$ is equal to one,
the Liouville's measure $d\mu = du_1...du_N$ is invariant under the action of $T$.
The inverse matrix $T^{-1}$ is also an integer matrix because $\Det  T=1$.
Therefore $T$ is an automorphism of  the unit hypercube  onto itself.
The conditions (\ref{mmatrix1}) on the eigenvalues of the matrix $T$ are  sufficient
to prove that the system represents  an Anosov C-system \cite{anosov} and therefore as such it
also represents a  Kolmogorov K-system \cite{kolmo,kolmo1,sinai3,rokhlin,rokhlin2}
with mixing  of  all orders and of nonzero entropy.

\begin{figure}
 \centering
\includegraphics[width=8cm]{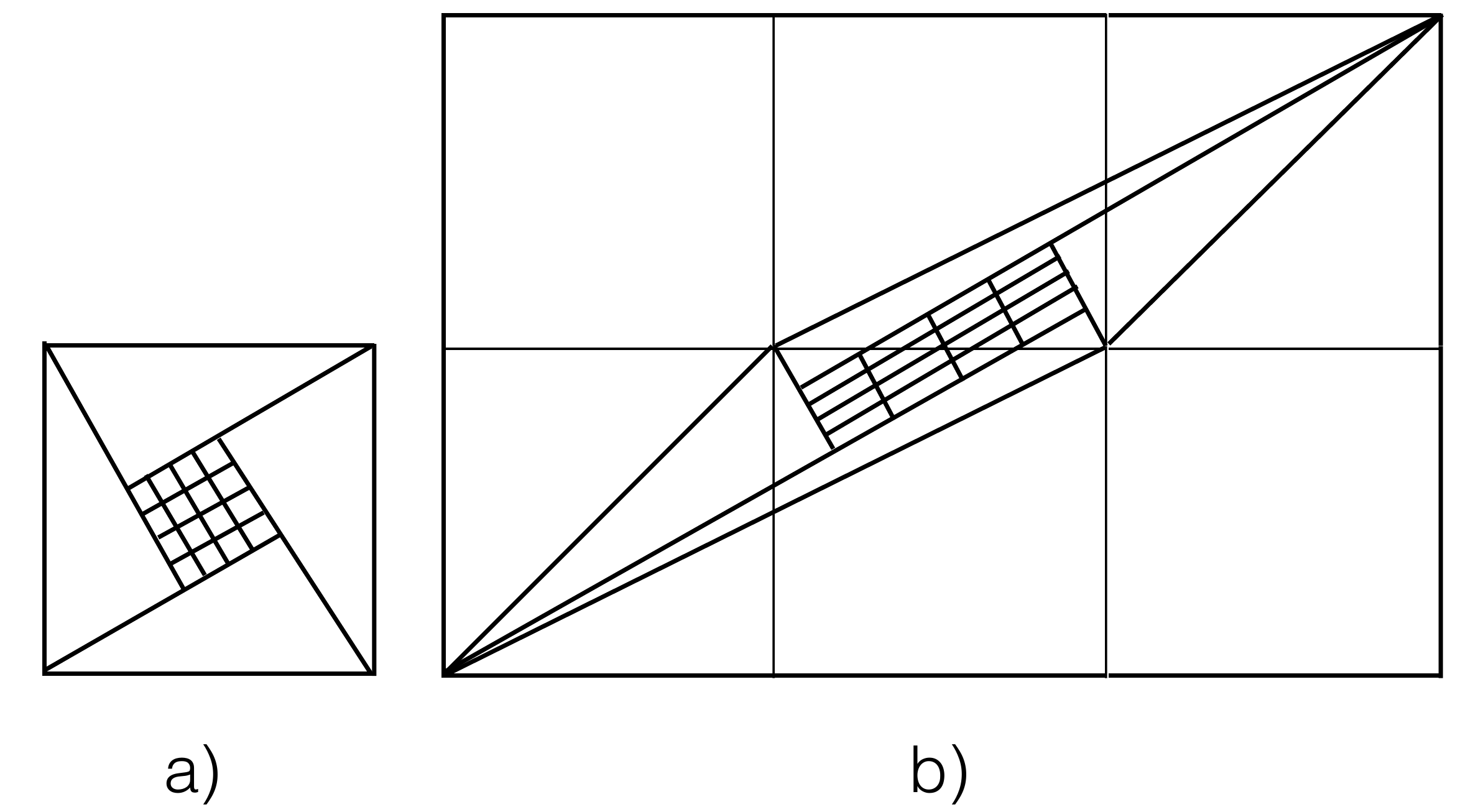}
\centering
\caption{The eigenvectors of the matrix $T$
 $\{ e_{\alpha}  \}$ and  $\{ e_{\beta}  \}$ define
two families of  parallel planes  $X_{\alpha} $ and $Y_{\beta}$,
which are invariant under the automorphisms  $T$.
The automorphism $T$ is contracting  the distances between points on the
planes belonging to the set $X_{\alpha}$ and expanding the
distances between points on the planes belonging to the set $Y_{\beta}$.
The a) depicts the parallel planes of the sets  $X_{\alpha} $ and $Y_{\beta}$
and b) depicts their positions after the action of the automorphism  $T$.}
\label{fig1}
\end{figure}

Let us recall how one can compute  the entropy $h(T)$ of the torus automorphisms (\ref{eq:rec}) 
$
u^{(n)} = T^n u^{(0)}
$
by using the eigenvalues of the  matrix $T$.  The eigenvalues of the matrix $T$ are divided
into two sets   $\{ \lambda_{\alpha}  \} $ and $\{  \lambda_{\beta }  \} $
with modulus smaller and larger than one:
\bea\label{eigenvalues}
&0 <  \vert \lambda_{\alpha} \vert   < 1~~~~   & \textrm{ for } \alpha=d{+}1,...,N \nn\\
&1 <  \vert \lambda_{\beta} \vert  < \infty & \textrm{ for } \beta=1,...,d . 
\eea
There exist two  hyperplanes $ X_{\alpha}$ and $ Y_{\beta} $,
which are  spanned by the
corresponding eigenvectors  $\{ e_{\alpha}  \}$ and  $\{ e_{\beta}  \}$ .
These invariant planes of the matrix $T$,  for which the eigenvalues are  inside and outside 
of the unit circle respectively,  define exponentially contracting  and expanding invariant subspaces. The phase trajectories are  expanding and contracting
under the transformation $T$ at an exponential rate (see Fig.\ref{fig1}).  The same is true for the inverse evolution that is defined by the matrix $T^{-1}$. For the inverse evolution the
contracting and expanding invariant spaces alternate their role.

\begin{figure}
 \centering
\includegraphics[width=3cm]{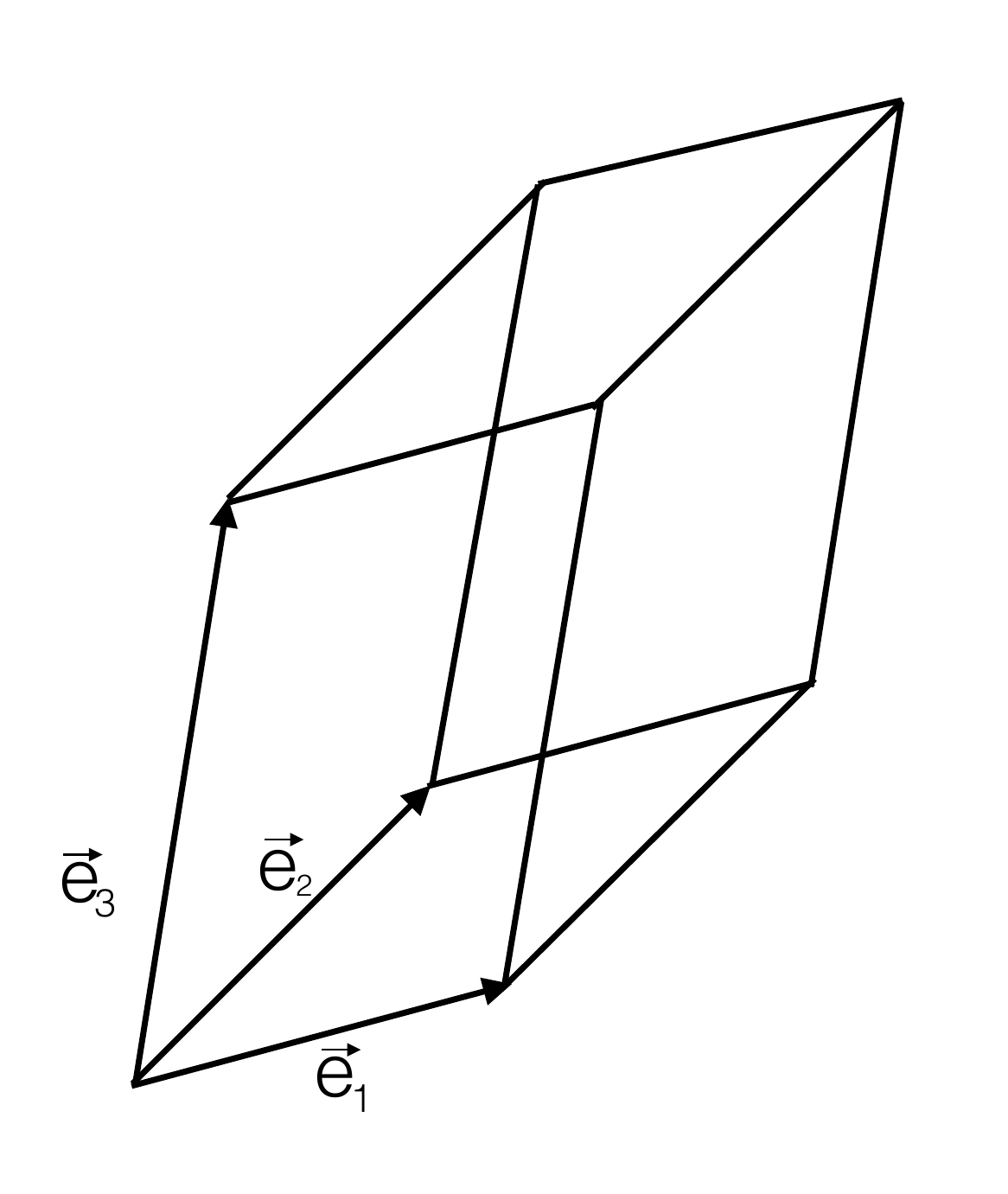}
\centering
\caption{A three-dimensional parallelogram ($d=3$) embedded into the expanding foliation $Y_{\beta}$. It is a simplex consisting of sites of lower-dimensional simplexes of dimension one, two and three.  The one-dimensional simplexes are the $\vec{e}_{i},~i=1,2,3$, the two-dimensional are constructed by the pair of vectors  $(\vec{e}_{i},\vec{e}_{j})~{i} \neq {j}$ and the three-dimensional one by $(\vec{e}_{1},\vec{e}_{2},\vec{e}_3)$. The extended entropy (\ref{extendedentropy}) is  $r_2(T) = \ln( \vert \lambda_1\vert) \ln (\vert\lambda_2 \lambda_3\vert) +  \ln (\vert \lambda_2 \vert) \ln ( \vert \lambda_3 \lambda_1\vert ) +\ln (\vert \lambda_3\vert ) \ln (\vert \lambda_1 \lambda_2\vert )$.}
\label{fig2}
\end{figure}

A convenient way to calculate the entropy $h(T)$  \cite{anosov,sinai3,rokhlin2,sinai4,gines,Savvidy:2015ida} is  to integrate over the whole phase space $M$ the logarithm of the volume expansion rate $\lambda(u)$ of an infinitesimal d-dimensional parallelogram (see Fig.\ref{fig1} and Fig.\ref{fig2} ) that is embedded into the expanding foliation $Y_{\beta}$ with its edges along the corresponding eigenvectors $e_{\beta}$, where 
$$
T  e_{\beta} =  \lambda_{\beta }  e_{\beta},~~~~\beta=1,...,d .
$$
For the automorphisms (\ref{eq:rec}) the coefficient $\lambda(u)$  does not depend on the phase
space coordinates $u$ and is equal to the product of eigenvalues
$
\{  \lambda_{\beta }  \} 
$  
with modulus  larger than one:
\be\label{volumeexp}
\lambda(u) = \prod^d_{\beta=1} \vert \lambda_{\beta} \vert .
\ee
Thus the  entropy of the automorphism  (\ref{eq:rec})  is  equal to the sum:
\be\label{entropyofA}
h(T) =\ln(\prod^d_{\beta=1} \vert \lambda_{\beta} \vert )= \sum^d_{\beta =1  } \ln \vert \lambda_{\beta} \vert.
\ee
It is clearly from this expression that its value depends on the eigenvalues of the evolution matrix $T$.  The expression (\ref{entropyofA}) for the entropy fulfils the Clausius conception of entropy to be extensive. Indeed, if the matrix $T$ has a block diagonal form with the matrices $T_1$ and $T_2$ on the diagonal, then  one can simply observe additivity of the entropy: 
\be\label{embedding}
h(T_1)+h(T_2) = \sum^{d_1}_{\beta_1 =1  } \ln \vert \lambda_{\beta_1} \vert +\sum^{d_2}_{\beta_2 =1  } \ln \vert \lambda_{\beta_2} \vert = \sum^d_{\beta =1  } \ln \vert \lambda_{\beta} \vert=  h(T)  .
\ee
The value of the Kolmogorov entropy $h(T)$ allows to quantitatively characterise the pseudorandom number generators \cite{yer1986a}. For a generator to pass a battery of statistical tests one should have a large enough entropy \cite{konstantin}.  What we observed when studying these generators was that even when the generators have equal and large enough entropies, nevertheless they generate a substantially different quality of random sequences. The question is how to differentiate these systems quantitatively when they have identical entropies? Is there an additional metrical invariant which will allow to define a refined classification of systems that have equal Kolmogorov entropies? 

In  \cite{yer1986a} it was conjectured that the generators based on hyperbolic systems (\ref{eq:rec}) of high dimensionality $N$  are advantageous compared to the low-dimensional ones. The argument, I quote, was:  "The advantage of a pseudorandom number generator, given by a dynamical system (\ref{eq:rec}), is that although the relaxation times 
\be\label{relaxation}
\tau_0 = 1/h(T)
\ee
may be made equal, the "quality" of mixing in the systems (\ref{eq:rec}) is higher owing to the fact, that in different directions the rate of instability is different and is proportional to the eigenvalues of the matrix $T$, which are quite arbitrary. Such "many-scale" mixing of directions ensures a slower growth of the Kolmogorov discrepancy $D_N(T)$."

It is true that the eigenvalues represent the metrical invariants of a system, but it is unclear how to compare  the eigenvalues of two independent hyperbolic systems.  In order to find a solution let us consider the low-dimensional subsystems  that are naturally embedded into a large one, like in (\ref{embedding}). In that case we have the low-dimensional expanding foliations that are embedded into the expanding foliation  $Y_{\beta}$  of a large system. The entropy of a low-dimensional subsystem is equal the logarithm of the expansion rate of the volumes of the low-dimensional parallelograms and are equal to the product of the corresponding eigenvalues (\ref{embedding}).  These quantities will characterise the entropies of the subsystems if a large system is factorised into a lower dimensional ones, that is when a matrix $T$ has block diagonal form  and we will have:
\be
h(T_1)+...+h(T_m)   = \sum^{d_1}_{\beta_1 =1  } \ln \vert \lambda_{\beta_1} \vert +...+\sum^{d_m}_{\beta_m =1  } \ln \vert \lambda_{\beta_m} \vert = h(T) .
\ee
Now one should understand how to proceed if the matrix $T$ does not have block diagonal form. In that case we still can consider lower-dimensional subsystems and calculate the corresponding entropies, but the reality is that  there are  many of them and none of them have  any "privileged" position within a large system. In that case it seems natural to consider all subsystems simultaneously. 

The d-dimensional parallelogram embedded into the expanding foliation $Y_{\beta}$ is a geometrical simplex of dimension $d $ consisting of sites of lower-dimensional simplexes: of dimension one, two and larger, up to the $d$-dimensional simplex (see Fig.\ref{fig2}).  The lowest-dimensional parallelograms will consist of $d$ one-dimensional simplexes (edges) parallel to the eigenvectors $e_{i},~i=1,...,d$. The logarithm of the expansion rate of the i-th parallelogram with give us the associated entropy $h_i = \ln \vert \lambda_i \vert $. The two-dimensional parallelograms will have two edges parallel to the eigenvectors  $(e_{i_1},e_{i_2})~{i_1} \neq {i_2}$ and its expansion rate will give the associated entropies $h_{i_1i_2} = \ln (\vert \lambda_{i_1}  \lambda_{i_2} \vert )$.  Considering the expansion rate of all possible distinct lower-dimensional parallelograms of increasing dimensionality we shall get the correspond "sub-entropies" 
\be
h_i = \ln \vert \lambda_i \vert,~~~ h_{i_1i_2} = \ln (\vert \lambda_{i_1}  \lambda_{i_2} \vert), ~~~.....~~~h_{i_1...i_m} = \ln (\vert \lambda_{i_1}....  \lambda_{i_m} \vert)~~...... 
\ee
The sum of all these entropies will give a quantity proportional to the standard entropy (\ref{entropyofA}), therefore it seems natural to construct a quantity that is quadratic in $h_{i_1...i_m}$, that is, to define an extended entropy as a sum $s(T) = h(T)  +r_2(T)...+r_d(T)$, where 
\bea\label{extendedentropy}
 r_2(T) &\sim& \sum_1 \ln \vert \lambda_{i_1} \vert  \cdot \ln \vert \lambda_{i_2}...\lambda_{i_d}  \vert +\nn \\
      &+& \sum_{2} \ln \vert \lambda_{i_1} \lambda_{i_2}  \vert  \cdot \ln \vert \lambda_{i_3}...\lambda_{i_d}  \vert +...\nn\\
      &+& \sum_{[{d\over 2}]} \ln \vert \lambda_{i_1} ...\lambda_{i_[{d\over 2}]}  \vert \cdot \ln \vert \lambda_{i_[{d\over 2}]+1}...\lambda_{i_d}  \vert 
= C(d) \sum^d_{i_1 \neq i_2}  \ln \vert \lambda_{i_1} \vert  \cdot \ln \vert \lambda_{i_2} \vert.  
   \eea 
and the sum $\sum_m$  runs over all non equal permutations of $i's$. We will simply define
\be\label{r2}
r_2(T) = \sum^d_{i_1 \neq i_2}   \ln \vert \lambda_{i_1} \vert  \cdot \ln \vert \lambda_{i_2} \vert.
\ee
The extended entropy $r_2(T)$ is a metrical invariant of the automorphisms $T$.  If the systems $S $ and $T$ are isomorphic $S = QT Q^{-1}$, where $Q$ is a unitary operator, then they have equal eigenvalues and  $r_2(S) =r_2(T)$. 
It follows from (\ref{r2}) that spectrally equivalent   dynamical systems $U_S = Q U_T Q^{-1} $ that have countable Lebesgue spectrum of the unitary  operators $U_S$ and $U_T$  and equal Kolmogorov entropies $h(S) =h(T)$ are not isomorphic when $r_2(S) \neq r_2(T)$. We shall define the high order entropies  as
\bea\label{r3-rd}
r_3(T) &=& \sum^d_{i_1 \neq i_2 \neq i_3}  \ln \vert \lambda_{i_1} \vert  \cdot \ln \vert \lambda_{i_2} \vert \cdot \ln \vert \lambda_{i_3} \vert \nn\\
....&&.......\nn\\
r_{d-1}(T) &=& \sum^d_{i_1 \neq i_2 \neq ...\neq i_{d-1}}  \ln \vert \lambda_{i_1} \vert  \cdot \ln \vert \lambda_{i_2} \vert \cdot...\cdot \ln \vert \lambda_{i_{d-1}} \vert \nn\\
\nn\\
r_d(T) &=&   \ln \vert \lambda_{1} \vert  \cdot \ln \vert \lambda_{2} \vert \cdot...\cdot \ln \vert \lambda_{d} \vert . 
\eea
The hyperbolic systems are not isomorphic even if all of these invariants are equal except one of them. In the next section we shall consider the relevant examples. 

But before considering examples it will be interesting to know if there is an equivalent definition of the extended entropies $r(T)$ that can be given in terms of the information theory \cite{kolmo,kolmo1,Shannon}. Let us in short recollect the construction of the Kolmogorov entropy $h(T)$.  If  $\alpha = \{A_i\}_{i \in I}$ ( $I$ is finite or countable) is a  measurable partition of the phase space $M$,~ 
$
\mu(M \setminus \bigcup_{i \in I} A_i)=0,~~~~\mu(  A_i \bigcap   A_j)=0, i \neq j~,
$
then the entropy of the partition $\alpha$ is
\be\label{partitionentropy}
h(\alpha) = - \sum_{i \in I} \mu(A_i) \ln \mu(A_i).
\ee
The {\it refinement partition} 
$
\alpha = \alpha_1 \vee \alpha_2  \vee ... \vee \alpha_k
$
of the 
collection of partitions $\alpha_1,..., \alpha_k$ 
is the intersection of all their composing sets  $A_i$:
$
\alpha = \big\{ \bigcap_{i \in I} A_i~ \vert ~A_i \in \alpha_i~ for ~all~ i  \big\}.
$
The entropy of the refinement 
$
\alpha \vee T \alpha \vee ...\vee T^{n-1} \alpha 
$
generated by iteration of the  automorphism $T$ 
is the limit:
\be
h(\alpha, T)= \lim_{n \rightarrow \infty} {h(\alpha \vee T \alpha \vee ...\vee T^{n-1} \alpha) \over n},~~~~
n=1,2,...,
\ee
and the entropy of the 
automorphism $T$ is a supremum taken over all finite measurable partitions $\{ \alpha \}$ of  $M$~\cite{kolmo,kolmo1,sinai3,rokhlin1,rokhlin,rokhlin2}: 
\be\label{supremum}
h(T) = \sup_{\{ \alpha \}} h(\alpha,T) .
\ee
In article \cite{Tsallis} Tsallis suggested an extension of the standard expression for the entropy
$
h=\sum^n_{i=1}p_i \ln {1\over p_i}
$
by using $q$ deformation of the logarithmic function: 
\be
h_q = {1- \sum^n_{i=1} p^{q}_i \over q-1} \equiv \sum^n_{i=1}p_i \ln_q {1\over p_i} ~ ,~~~~~~~\sum^n_{i=1}p_i =1.
\ee
The standard entropy is recovered when the parameter $q \in R$ tends to one $\lim_{q \rightarrow 1} h_q =h$. 
By using the $q$ deformation of the logarithm one can suggest the generalisation of the entropy of a partition $\alpha$ in (\ref{partitionentropy}) to be promoted to the  expression
\be
h_q(\alpha) =    \sum_{i \in I} \mu(A_i) \ln_q {1 \over \mu(A_i) }~.
\ee
The entropy of  the partition $\alpha$ with respect to the automorphisms $T$ is the limit
\be
h_q(\alpha, T)= \lim_{n \rightarrow \infty} {h(\alpha \vee T \alpha \vee ...\vee T^{n-1} \alpha) \over n},~~~~
n=1,2,...
\ee
and the generalised entropy of the automorphism $T$ we shall define as a supremum: 
\be\label{supremum1}
h_q(T) = \sup_{\{ \alpha \}} h(\alpha,T).
\ee
The extended entropy $h_q(T)$  is a metrical invariant. Indeed, it is formulated in terms of a measure that is invariant under automorphism transformations  $\mu( A)= \mu(T A)$. If two systems $S $ and $T$ are isomorphic $S = QT Q^{-1}$, where $Q$ is a unitary map between phase spaces $K$ and $M$, then to each measurable partition $\alpha$ of  $M$ correspond a particular measurable partition $\beta$ of $K$,  therefore isomorphic systems will have equal values of  $h_q$. The new feature of the $h_q(T)$ consists in the fact that for two independent automorphisms $T_1$ and $T_2$ we will have the expression which is quadratic in subsystems entropies\footnote{The $q$ deformation of the logarithm has the following property: $\ln_q(x y) = \ln_q x  + \ln_q y +(1-q) \ln_q x  \ln_q y  $.}:
\be\label{alternativedef}
h_q(T) = h_q(T_1) + h_q(T_2) + (1-q)~ h_q(T_1)h_q(T_2).
\ee 
The last term is quadratic in the entropies of the subsystems and has a structure similar to the terms in the extended entropy $r_2(T)$. We were unable to derive our formulas (\ref{r2}),  (\ref{r3-rd}) from the alternative definition (\ref{supremum1}).  But it is appealing to conjecture that it is 
the $q$ deformation logarithm  of the volume expansion rate $h_q(T) = \ln_q \lambda(u)$.  That leads to the expression  
\bea
h_q(T) = \ln_q(\prod^d_{\beta=1} \vert \lambda_{\beta} \vert ) &=&\sum^{d}_{i_1=1}   \ln_q \vert \lambda_{i_1} \vert  +\nn\\
&+&(1-q) \sum^d_{i_1 \neq i_2} \ln_q \vert \lambda_{i_1}  \vert  \cdot  \ln_q \vert  \lambda_{i_2}  \vert \nn\\
&+&(1-q)^2  \sum^d_{i_1 \neq i_2 \neq i_3 }  \ln_q \vert \lambda_{i_1}  \vert  \cdot  \ln_q \vert  \lambda_{i_2} \vert \cdot  \ln_q \vert  \lambda_{i_3}  \vert+...\nn\\
...&+&(1-q)^{d-1}   \ln_q \vert \lambda_{1} \vert  \cdot \ln_q \vert \lambda_{2} \vert \cdot...\cdot \ln_q \vert \lambda_{d} \vert , 
\eea
which is identical in structure with the linear sum $s(T) = h(T) +r_2(T)+....+r_{d}(T)$ (\ref{r2}),  (\ref{r3-rd}). It is a challenging problem to prove this relation.
 
\section{\it Calculating Extended Entropy}

We shall consider the  MIXMAX matrix $T$ of the form \cite{yer1986a,konstantin}:
\be
\label{eqmatrix}
T(N,s) =  
   \begin{pmatrix}  
      1 & 1 & 1 & 1 & ... &1& 1 \\
      1 & 2 & 1 & 1 & ... &1& 1 \\
      1 & 3+s & 2 & 1 & ... &1& 1 \\
      1 & 4 & 3 & 2 &   ... &1& 1 \\
      &&&...&&&\\
      1 & N & N-1 &  N-2 & ... & 3 & 2
   \end{pmatrix}
\ee
The matrix is constructed so that its entries are increasing together with the size $N$ of the matrix, and we have a family of matrices which are parametrised by the integers $N$ and $s$ and have 
determinant equal to one. It is defined recursively, 
since the matrix of size $N+1$ contains in it the matrix of the size $N$.  
In order to generate pseudorandom vectors  $u^{(n)} = T^n u^{(0)}$, one should 
 choose  the initial vector $u^{(0)}=(u^{(0)}_{1},...,u^{(0)}_{m})$, called the "seed", 
 with at least one non-zero component to avoid fixed point of $T$, which is at the origin.
 
\begin{figure}
\begin{center}
\includegraphics[width=5cm]{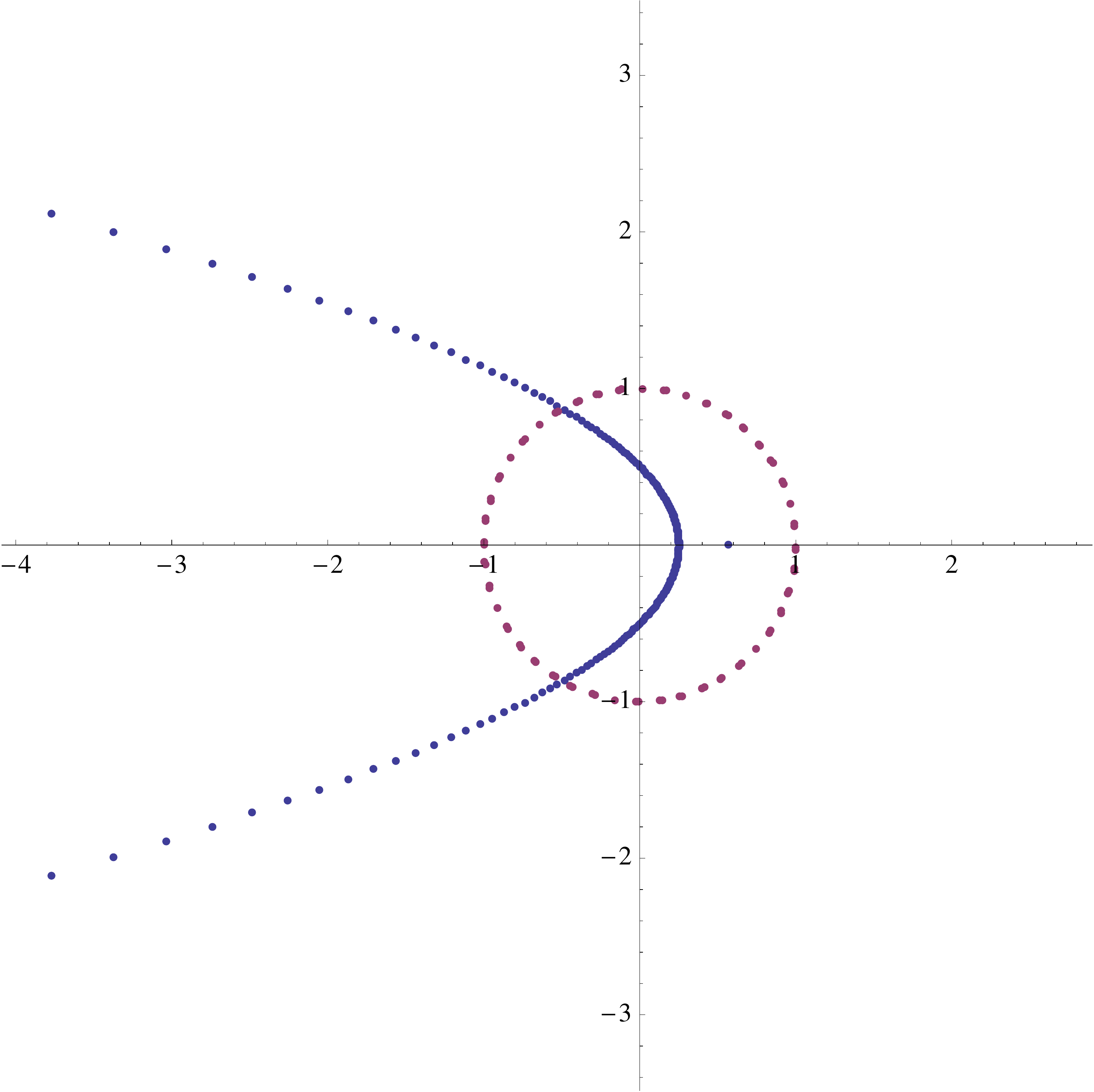}~~~~~~~~~~
\includegraphics[width=5cm]{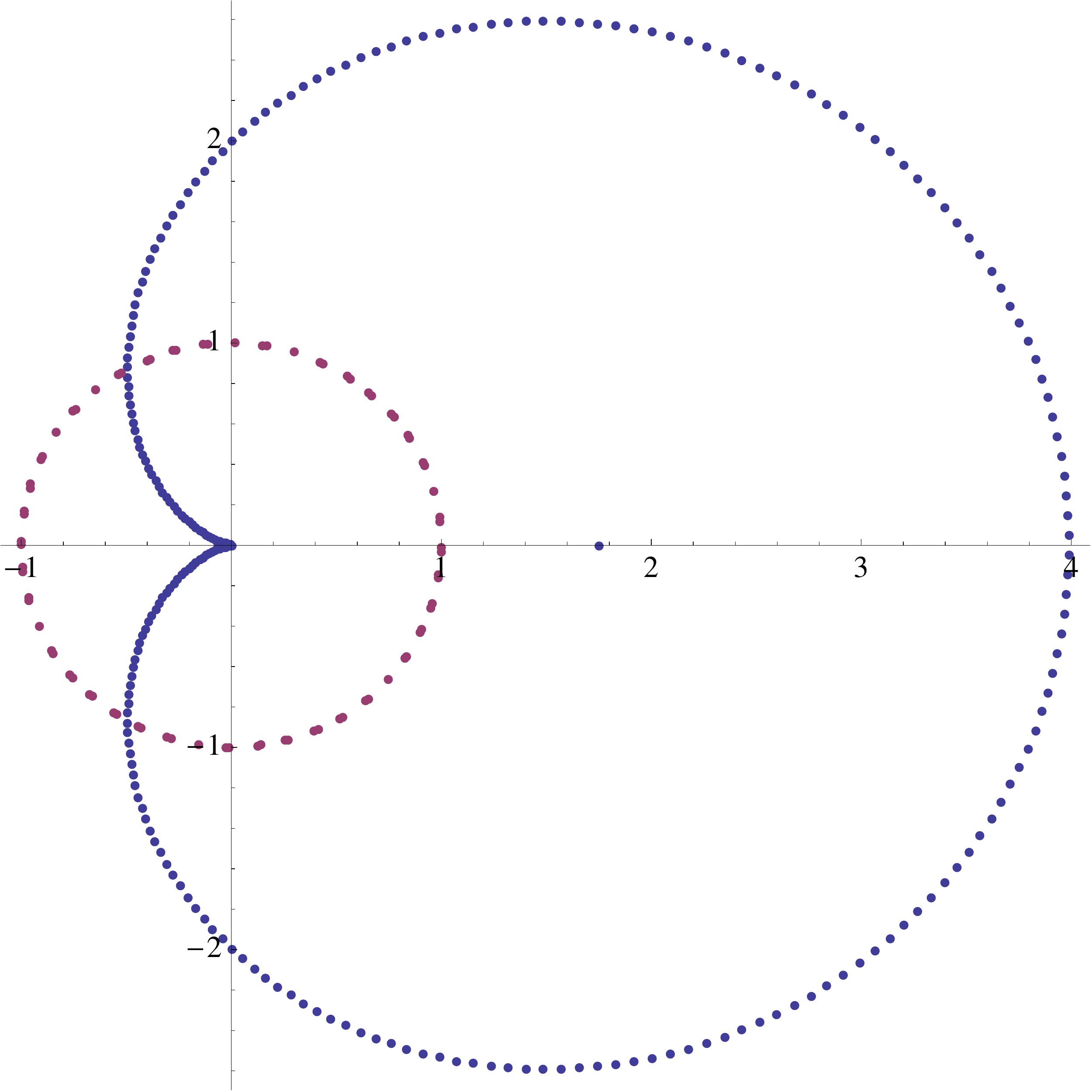}
\caption{
On the left is the distribution of the eigenvalues of the matrix $T$ in (\ref{eqmatrix})
and on the right of its  inverse matrix $T^{-1}$.
The unit circle is depicted to separate  the eigenvalues inside and outside
the circle in  accordance with  the formula (\ref{eigenvalues}). The eigenvalues 
of the matrices  $T$ and $T^{-1}$ lie on a parabola and on a cardioid.
A cardioid is the inverse curve of a parabola with its focus at the centre of inversion.
}
\label{fig3}
\end{center}
\end{figure}

 The eigenvalues of the matrices  $T$ and $T^{-1}$ (\ref{eqmatrix}) are complex  valued $ \lambda_i = r(\phi_i) \exp(i\phi_i)$ and lie on the parabola and cardioid correspondingly when $N\rightarrow \infty$\cite{konstantin,Savvidy:2015jva}: 
\be\label{eq:curve}
r(\phi) = {1 \over 4 \cos^2(\phi/2)},~~~~~~~~~r(\phi) = 4 \cos^2(\phi/2).
\ee
For finite $N$, the formula for the eigenvalues is
\be
\lambda_j = \frac{1}{4 \cos^2(j\pi/2N)} ~~ \exp(i\, \pi j/N) ~~ \textrm{for}~ j=-N/2..N/2 ~,
\label{eq:evs}
\ee
where the complex conjugate eigenvalues correspond to $\pm j$. 
The eigenvalues are widely dispersed and the spectrum is indeed  "multi-scale" as required  \cite{yer1986a} and is shown on  Fig.\ref{fig3}. From the above analytical expression for eigenvalues it follows that the eigenvalues satisfying the condition $ 0 <  \vert \lambda_{\alpha} \vert   < 1$  are in the range $-2\pi/3 <  \phi < 2\pi/3$ and the ones satisfying the condition $ 1 <  \vert \lambda_{\beta} \vert $ are in the interval $2\pi/3 <  \phi < 4\pi/3$.  

The entropy of the system $h(T)$ can be calculated for
large values of $N$  as an integral over eigenvalues (\ref{eq:curve}) \cite{Savvidy:2015jva}:
\be\label{linear}
h(T)=   \sum_{-2\pi/3 <  \phi_i < 2\pi/3} \ln (4\cos^2(\phi_i/2)~ \approx  ~N \int^{2\pi/3}_{-2\pi/3} \ln (4\cos^2(\phi/2){d\phi \over 2\pi} \approx
4.06 \Big({N\over 2\pi}\Big)
\ee
and grows linearly with the dimension $N$ of the operator $T$.
Let us now calculate the extended entropy (\ref{extendedentropy}) of the system  $T$ (\ref{eqmatrix}):\bea\label{Nextendedentropy}
& r_2(T) \geq ~ \sum_{-2\pi/3 <  \phi_m < 2\pi/3} ~\Big(~\sum_{-2\pi/3 <  \phi_i < \phi_m} \ln (4\cos^2(\phi_i/2)~\Big)~\Big( \sum_{\phi_m <  \phi_i < 2\pi/3} \ln (4\cos^2(\phi_i/2)\Big) \approx \nn\\
&\nn\\
& \approx N^3  \int^{2\pi/3}_{-2\pi/3}  ~ {d\phi \over 2\pi} ~
\Big(  \int^{\phi}_{-2\pi/3}   \ln (4\cos^2(\chi/2){d\chi \over 2\pi}    ~\Big)~       
     \Big( \int^{2\pi/3}_{\phi}  \ln (4\cos^2(\omega/2){d\omega \over 2\pi} ~\Big) \approx \nn\\
&      \nn\\
   &   \approx   9.138 \Big({N \over 2\pi} \Big)^3.
 \eea 
The extended entropy is growing  as  the cube of the matrix dimension $N$. The inequality appears due to the fact that not all of the permutations are taken into account in the integral representation (\ref{Nextendedentropy}).    

In Table \ref{tbl:largeN} we present the entropies of the operator $T(N,s)$ for the matrices of a large matrix size $N$. The third column represents the values of the Kolmogorov entropy $h(T)$ and the forth one the extended entropies $r(T)$, which are fairly  large. This reflects the fact that the eigenvalues are large and are well distributed.

\begin{table}[htbp]
   \centering
   \begin{tabular}{@{} lcccrcl @{}} 
         \toprule
      Size & Magic & Entropy &Extended Entropy& Log of the period $q$ \\
      N    & $s$ & h(T)  &$r_2(T)$ & $\approx \log_{10} (q)$  \\ 
      \midrule
        256   & -1 & 164.5& 618061  & 4682 \\      
        7307   & 0 & 4676.5& 1.4 $10^{10}$ & 134158 \\
	20693 & 0 & 13243.5& 3.3 $10^{11}$ & 379963 \\
        25087 & 0 & 16055.7 & 5.8 $10^{11}$& 460649 \\
 	28883 & 1 & 18485.1 & 8.9 $10^{11}$& 530355 \\
	40045 & -3 & 25628.8 & 2.4 $10^{12}$& 735321 \\
	44851 & -3 & 28704.6 &3.3 $10^{12}$& 823572 \\
      \bottomrule
   \end{tabular}
   \caption{Table of parameters of the operator $T(N,s)$ for large matrix size $N$. The third column is the value of the Kolmogorov entropy and the forth one of the extended entropy. The $\log_{10} q$ is logarithm of the period q. }
   \label{tbl:largeN}
\end{table}

Let us also consider some of the other popular generators. In the case of RCARRY \cite{RCARRY}, which is a slight modification of a Fibonacci-like recurrence modulo $2^{24}$, its failure was related to the weak mixing properties of its underlying matrix, because the entropy is of order $h \approx 0.32$. Its extended entropy is  $r \approx 1.25$. In case of the skipping procedure  these parameters are increasing by number of skipping, which is about hundred, but still remain small. These invariants remain small also for the  Mersenne Twister (MT) generators \cite{MT}.

\section{\it Conclusion}

The efficient implementation of the  C-system MIXMAX generators (\ref{eq:rec}) for Monte-Carlo simulations can be found in the articles \cite{konstantin,hepforge}.  
The MIXMAX generators demonstrated excellent statistical properties, high performance and superior high quality output and became a multidisciplinary usable product. The main characteristics of the generators are: a) MIXMAX is an original and genuine 64-bit generator, is one of the fastest generators producing 64-bit pseudorandom number in approximately 4 nanoseconds,
b) has very large Kolmogorov entropy of 0.9 per/bit,
c) long periods of order of $10^{120}$ - $10^{5000}$,
d) a new skipping algorithm generates seeds and guarantees that streams are not overlapping. The MIXMAX generators were integrated into the concurrent and distributed MC toolkit Geant4 \cite{Geant4}, the foundation library CLHEP \cite{CLHEP} and data analysis framework ROOT \cite{root}. These software tools have wide applications in High Energy Physics at CERN, in CMS experiment \cite{CMSrunII,CMS}, at SLAC, FNAL and KEK National Laboratories and are part of the CERN's active Technology Transfer policy. The generator is available in the PYTHIA event generator \cite{PYTHIA}. The MIXMAX code can be downloaded from the GSL-GNU Scientific Library \cite{GSL}.

\section{\it Acknowledgement }
I would like to thank   K. Savvidy  for stimulating discussions and long-lasting collaboration.

 \section{\it Appendix A}

The automorphism $T$ of the phase spaces $M$ is metrically isomorphic to the automorphism $S$ on $K$ if there exists the isomorphism $Q$ of the spaces  $M$ and $K$ such that $S = Q^{-1}T Q$. The spectral properties of the automorphism $T$ on the Hilbert space $L_2(M)$ of functions $\{f(u) ; u \in  M\}$ is defined as the spectrum of the corresponding conjugate unitary operator $U_{T}$  \cite{Koopman}:
\be\label{geoflow}
U_{T} f(u ) = f(T u), ~~~u \in  M.
\ee
The automorphisms $S$ and $T$ are considered spectrally equivalent if the unitary operators $U_{S}$ and $U_{T}$ have identical spectrum, that is $U_S = Q^{-1}U_T Q $.

 \section{\it Appendix B}

The C-condition  was  formulated by Anosov in  \cite{anosov}.
A {\it cascade}  on the m-dimensional compact phase space $W^{m}$ is induced by 
the diffeomorphisms $T: W^m \rightarrow W^m$.  The iterations
are defined by  a repeated action of the operator  
$\{ T^n, -\infty < n < +\infty  \}$, where $n$ is an  integer number.  
The tangent space at 
the point $w \in W^m$ is  denoted by $R^m_{w}$ and the 
tangent vector bundle by $\CR(W^m)$. 
The diffeomorphism $\{T^n\}$ induces the mapping
of the tangent spaces $\tilde{T}^{n}: R^m_w \rightarrow R^m_{ T^{n}w}$.
\begin{figure}
  \centering
\includegraphics[width=6cm]{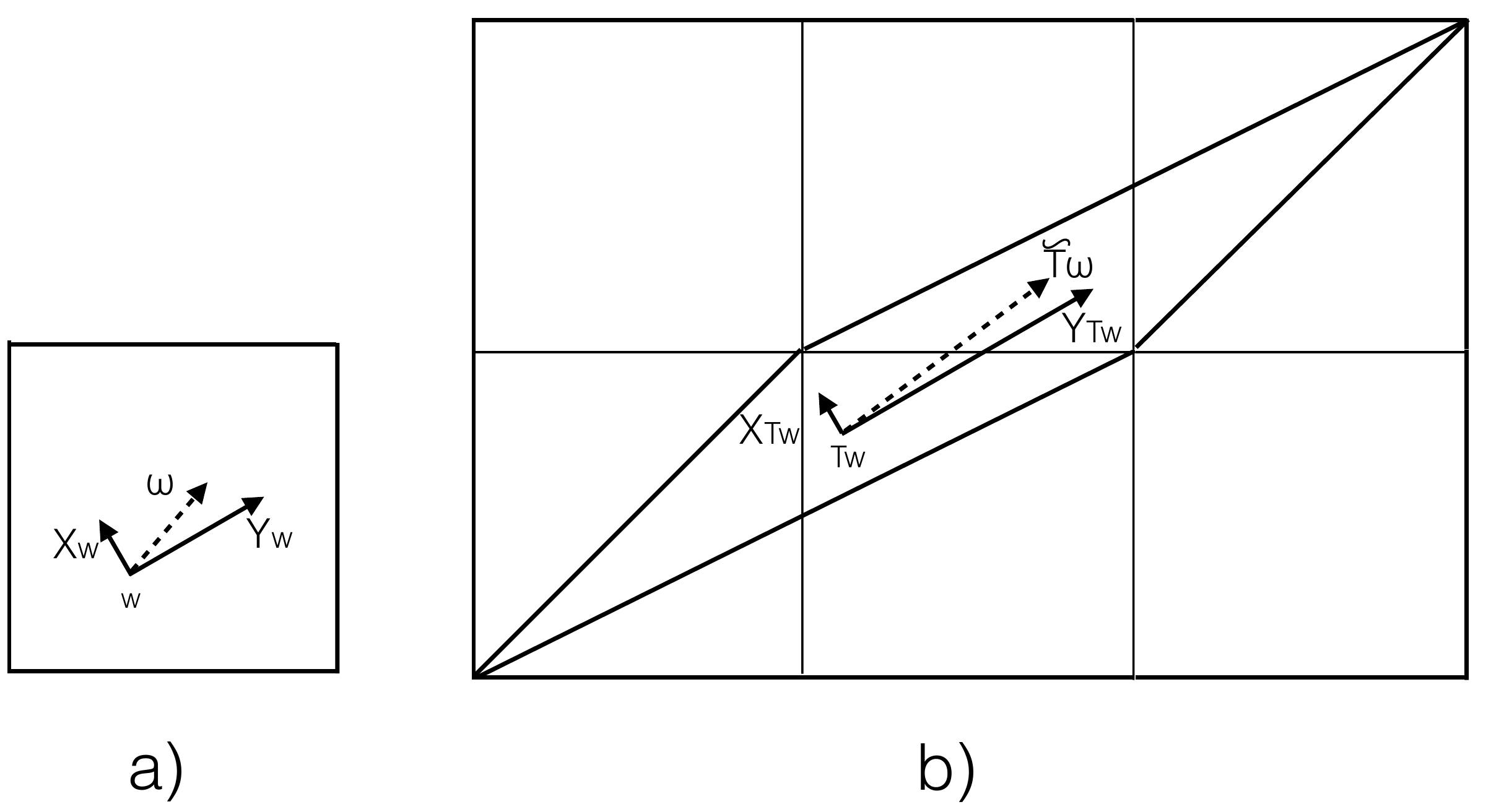}
\centering
\caption{Any tangent vector $\omega \in R_{w}$ at 
the point $w \in \mathbb{E}^N$  is decomposable  
into the sum $R_{w}= X_{w} \bigoplus Y_{w}$
where the spaces $X_w$ and $Y_w$ are defined by the corresponding 
eigenvectors of the matrix T (\ref{eq:rec}) .  
The automorphisms $T$ induces the mapping
of the tangent spaces $\tilde{T} X_{w} = X_{T w}, ~
\tilde{T} Y_{w} = Y_{T w}$.
It is contracting the distances on  $ X_{w}  $ and expanding the 
distances on  $ Y_{w} $.} 
\label{fig5}
\end{figure}
The C-condition requires that the tangent space $R^{m}_{w}$ at each point $w$
of the m-dimensional phase space $W^{m}$ of the dynamical system $\{T^{n}\}$ 
should be decomposable  into a 
direct sum of the two linear spaces  $X^{k}_{w}$ and $Y^{l}_{w}$ with the following 
properties \cite{anosov}:
\bea\label{ccondition}
C1.&R^{m}_{w}= X^{k}_{w} \bigoplus Y^{l}_{w} ~~\\
&The~ dynamical ~system~ \{T^{n}\}~ is ~such ~that: \nn\\
C2.&~~a) \vert \tilde{T}^{n} \xi  \vert  \leq ~ a \vert   \xi \vert e^{-c n}~ for ~n \geq 0 ; ~
\vert \tilde{T}^{n} \xi  \vert  \geq~ b \vert \xi \vert e^{-c n} ~for~ n \leq 0,~~~\xi \in  X^{k}_{w}, \nn\\
&b) \vert \tilde{T}^{n} \eta  \vert  \geq~ b \vert \eta \vert e^{c n} ~~for~ n \geq 0;~
\vert \tilde{T}^{n} \eta  \vert  \leq ~ a \vert   \eta \vert e^{c n}~ for~ n \leq 0,~~~\eta \in Y^{l}_{w},\nn
\eea
where the constants a,b and c are positive and are the same for all $w \in W^m$ and all 
$\xi \in  X^{k}_{w}$, $\eta \in Y^{l}_{w}$.
The length $\vert ...\vert$ of the tangent vectors   $\xi $ and $  \eta $  
is defined by the Riemannian metric $ds$ on $W^m$.
The linear spaces $X^{k}_{w}$ and $Y^{l}_{w}$ are invariant  with respect to 
the derivative  mapping  $\tilde{T}^{n} X^{k}_{w} = X^{k}_{T^n w}, ~
\tilde{T}^{n} Y^{l}_{w} = Y^{l}_{T^n w}$ and represent the {\it contracting and expanding 
linear spaces} (see Fig.\ref{fig5}).
The C-condition describes the behaviour of all trajectories $\tilde{T}^n \omega$ 
on the tangent vector bundle  $\omega \in R^{m}_{w}$. 
Anosov proved that the vector spaces  $X^{k}_{w}$ and $Y^{l}_{w}$ are continuous 
functions of the coordinate $w$ and that they are the target vector spaces  to 
the foliations $\Sigma^k$ and $\Sigma^l$ which are  the {\it  surfaces transversal to 
the trajectories}  $T^n w$ on $W^m$ (see Fig. \ref{fig5}). 
\begin{figure}
 \centering
\includegraphics[width=6cm]{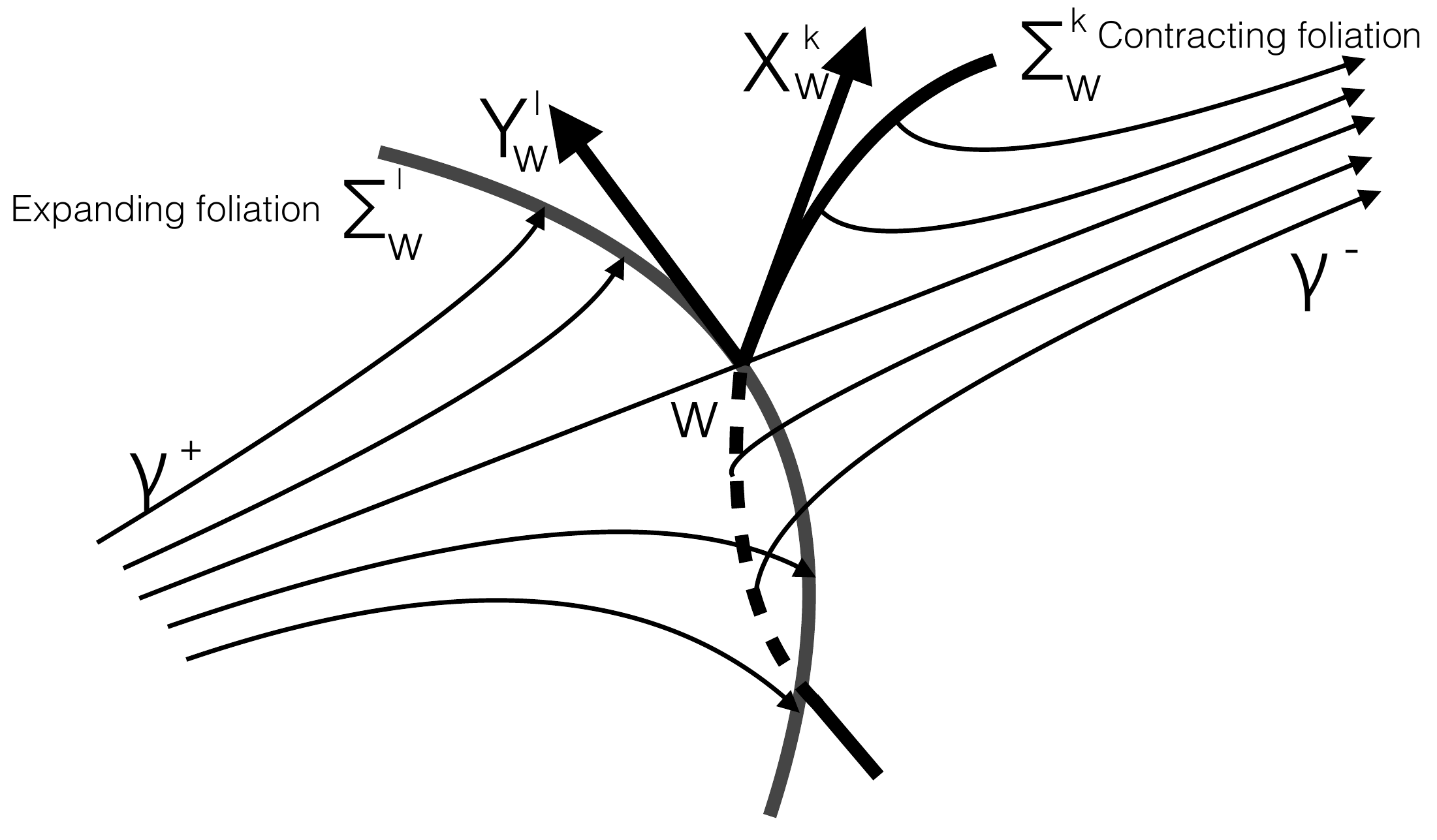}
\centering
\caption{ At each point $w$ of the C-system the tangent space $R^{m}_{w}$  
is  decomposable  into a direct sum of two linear spaces  $Y^l_{w} $ and $X^k_{w} $. 
The expanding and contracting geodesic flows
are $\gamma^+$ and $\gamma^-$. The expanding and 
contracting invariant  foliations  $\Sigma^l_{w} $ and $\Sigma^k_{w} $ 
{\it are transversal to the geodesic flows} and their corresponding tangent spaces 
are  $Y^l_{w} $ and $X^k_{w} $.  } 
\label{fig5}
\end{figure}

\section{\it Appendix C}

\begin{figure}
 \centering
\includegraphics[width=7cm]{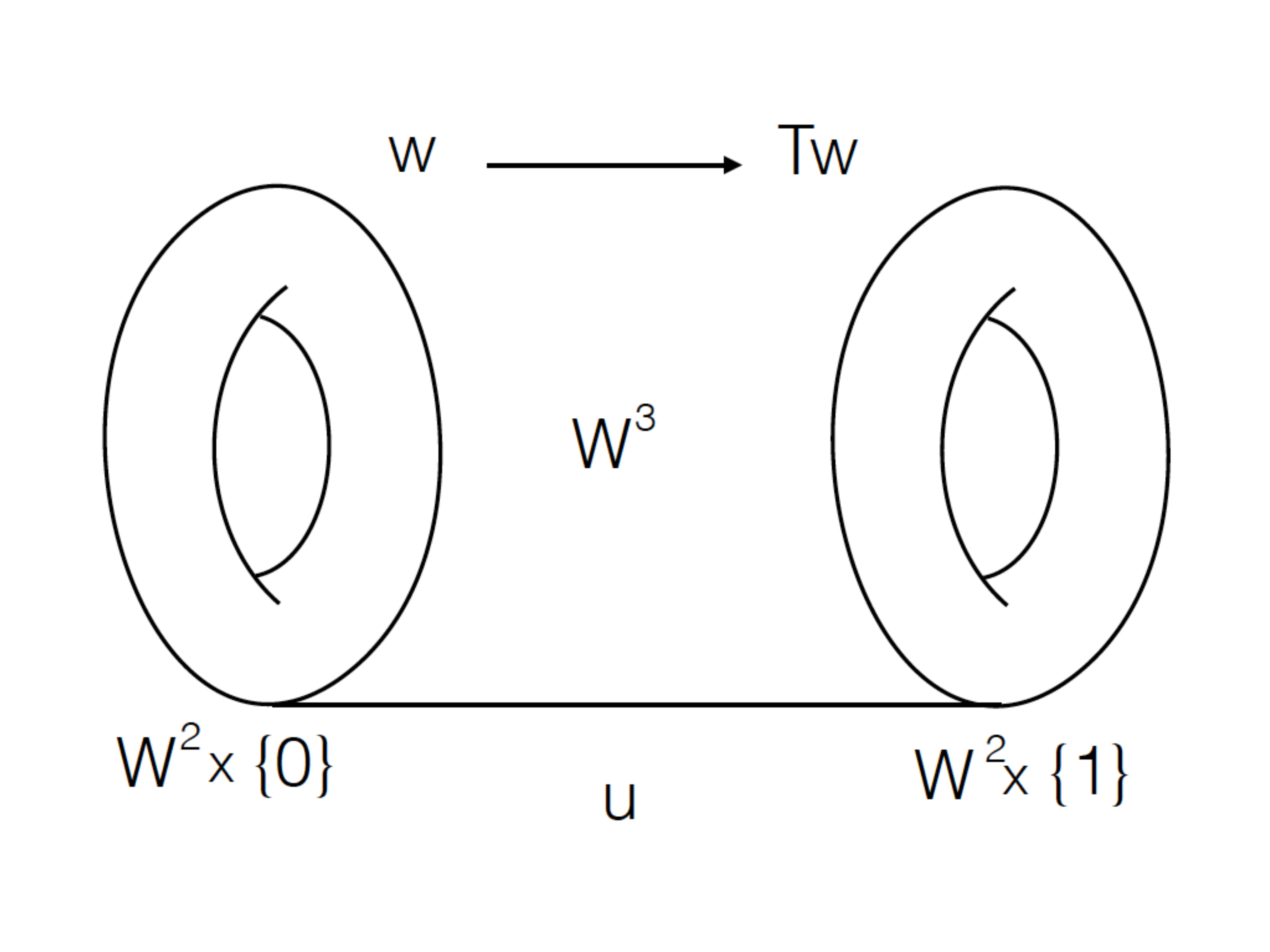}
\caption{The  identification of the fibres $W^2 \times \{0\}$ and  $W^2 \times \{1\}$ is made by the 
formula $(w,1) \equiv  (Tw,0)$ and the resulting manifold is a closed cylinder $W^2 \times [0,1]$, where $[0,1]= \{  u~ \vert ~ 0 \leq u \leq 1 \}$.
The resulting compact manifold $W^{3}$ has a bundle structure 
with the base $S^1$ and fibres $W^2 \times \{u\}$. The manifold $W^{3}$ 
has the local coordinates $\tilde{w}=(w^1,w^2,u)$ . 
} 
\label{fig6}
\end{figure}

 Any C-cascade on a torus can be embedded into a certain  C-flow \cite{anosov}. 
Let us consider a C-cascade on a torus $W^m$
(defined in section two) and increase its dimension m  by one unit
constructing a cylinder $W^m \times [0,1]$, where $[0,1]= \{  u~ \vert ~0 \leq u \leq 1 \}$,
and identifying $W^m \times \{0\}$ with $W^m \times \{1\}$ by the formula:
$
(w,1) \equiv  (Tw,0).
$
Here T is diffeomorphism (\ref{eq:rec}).
The resulting compact Riemannian manifold $W^{m+1}$ has a bundle structure
with the base $S^1$ and fibres of the type $W^m$ (see Fig.\ref{fig6}). The manifold $W^{m+1}$
has the local coordinates $\tilde{w}=(w^1,...,w^m,u)$ . The C-flow $T^t$ on the manifold $W^{m+1}$ is defined by the equations \cite{anosov}
\be\label{velo}
{d  w^1 \over d t}=0~, ....,~ {d  w^m  \over d t} = 0,~ {d  u \over d t}=1.
\ee
For this flow the tangent space $R^{m+1}_{\tilde{w}}$ is
a direct sum of three subspaces:
 \be
R^{m+1}_{\tilde{w}} = X^k_{\tilde{w}}  \oplus  Y^l_{\tilde{w}}  \oplus  Z_{\tilde{w}}.
\ee
The linear space  $X^k_{\tilde{w}} $ is tangent to the fibre  $W^m \times u$ and is parallel
to the eigenvectors corresponding to the eigenvalues which are lying inside the unit circle  $0 <  \vert \lambda_{\alpha} \vert   < 1$ and $Y^l_{\tilde{w}} $ is tangent to the fibre  $W^m \times u$ and is parallel to the eigenvectors  corresponding to the eigenvalues
which are lying outside of the unit circle $1 <\vert \lambda_{\beta}\vert$. $ Z_{\tilde{w}}$ is collinear to the phase space velocity (\ref{velo}). Under the derivative mapping of the (\ref{velo}) the vectors from $X^k_{\tilde{w}} $ and  $Y^l_{\tilde{w}} $ are contracting  and
expanding:
\be
\vert \tilde{T}^{t} e_{\alpha} \vert =  \lambda_{\alpha}^{ t}~ \vert e_{\alpha} \vert,~~~~
\vert \tilde{T}^{t} e_{\beta} \vert =  \lambda_{\beta}^{ t}~ \vert e_{\beta} \vert.
\ee
This identification of contracting and expanding spaces proves  that (\ref{velo})
defines a C-flow.

We are interested now to define a {\it geodesic flow} on the same Riemannian manifold $W^{m+1}$. 
The geodesic flow on $W^{m+1}$ was not studded in \cite{anosov} and has dynamics which is different  from (\ref{velo}). The equations for the geodesic flow on $W^{m+1}$
\be
{d^2 \tilde{w}^{\mu} \over d t^2} +\Gamma^{\mu}_{\nu \rho}
{d  \tilde{w}^{\nu} \over d t } {d  \tilde{w}^{\rho} \over d t } =0
\ee
are different from the flow equations defined by the equations (\ref{velo}).  If all sectional
curvatures on $W^{m+1}$ are negative then geodesic flow defines a C-flow \cite{anosov}. For simplicity we shall consider a  two-dimensional case $m=2$ and the matrix
\be\label{arno}
T= \begin{pmatrix}
1&1\\
1&2\\
 \end{pmatrix} .
 \ee
The metric on the $W^{m+1}$ is defined as:
\be\label{metric}
ds^2 = e^{2u} [\lambda_1 d w^1 + (1-\lambda_1) d w^2]^2 +
e^{2u} [\lambda_2 d w^1 + (1-\lambda_2) d w^2]^2 +du^2 =\nn\\
 g_{\mu\nu} d\tilde{w}^{\mu} d\tilde{w}^{\nu},
\ee
where $0 < \lambda_2  < 1 <  \lambda_1$ are eigenvalues of the matrix (\ref{arno}) and
fulfil the relations $\lambda_1 \lambda_2 =1,\lambda_1+ \lambda_2=3$.
The metric  is invariant  under the automorphism $T$
\be\label{trans}
w^1 = 2 w^{'1} - w^{'2},~~~ w^2 = -w^{'1}_1 + w^{'2}, ~~~u  = u^{'}-1
\ee
and is therefore consistent with the identification $(w,1) \equiv  (Tw,0)$. The metric
tensor has the form
\be\label{metric}
 g_{\mu\nu}(u)= \begin{pmatrix}
\lambda_1^{2 + 2 u} + \lambda_2^{2 + 2 u} & (1 - \lambda_1) \lambda_1^{1 + 2 u} + (1 - \lambda_2) \lambda_2^{1 + 2 u}& 0 \\
(1 - \lambda_1)\lambda_1^{1 + 2 u} + (1 - \lambda_2) \lambda_2^{1 + 2 u}&(1 - \lambda_1)^2 \lambda_1^{2 u} + (1 - \lambda_2)^2 \lambda_2^{2 u}&0\\
0&0&1\\
 \end{pmatrix}
 \ee
and the corresponding geodesic equations take the following form:
 \bea\label{geodesicflow}
&  \ddot{w}^1 + 2{(\lambda_1-1) \ln\lambda_1 \over \lambda_1+1} \dot{w^1}\dot{u}
-4 {(\lambda_1-1) \ln\lambda_1 \over \lambda_1+1} \dot{w^2}\dot{u}  =0
\nn\\
&  \ddot{w}^2 - 2{(\lambda_1-1) \ln\lambda_1 \over \lambda_1+1} \dot{w^2}\dot{u}
  - 4 {(\lambda_1-1) \ln\lambda_1 \over \lambda_1+1} \dot{w^1}\dot{u} =0\\
& \ddot{u} + {(1-\lambda^{4u+4}_1) \ln\lambda_1 \over \lambda^{2u+2}_1} \dot{w^1}\dot{w^1}
 + 2{(1+\lambda^{4u+3}_1)(\lambda_1-1) \ln\lambda_1 \over \lambda^{2u+2}_1} \dot{w^1}\dot{w^2}+\nn\\
& +  {(1-\lambda^{4u+2}_1)(\lambda_1-1)^2 \ln\lambda_1 \over \lambda^{2u+2}_1} \dot{w^2}\dot{w^2}=0.\nn
\eea
One can become convinced that these equations are invariant under the transformation (\ref{trans}).
 In order to study a stability of the geodesic flow one has to compute the sectional curvatures.
We shall choose the orthogonal frame in the directions of the linear spaces  $X^1_{\tilde{w}} , Y^1_{\tilde{w}} $   and $ Z_{\tilde{w}}$. The corresponding vectors are:
$
e_1 = (\lambda_1-1, \lambda_1,0),~e_2=(\lambda_2 -1, \lambda_2,0),~ e_3=(0,0,1)
$
and in the metric (\ref{metric}) they have the lengths:
\be
\vert e_1  \vert^2= (\lambda_1 - \lambda_2)^2 \lambda_2^{2 u},~~~~
\vert e_2  \vert^2= (\lambda_1 - \lambda_2)^2 \lambda_1^{2 u},~~~~
\vert e_3 \vert^2 = 1.
\ee
The corresponding sectional curvatures are:
\bea
K_{12} = {R_{\mu\nu\lambda\rho} e^{\mu}_1 v^{\nu}_2  e^{\lambda}_1 e^{\rho}_2 \over
\vert e_1 \wedge e_2 \vert^2}=   \ln^2 \lambda_1 >0
\nn\\
K_{13} = {R_{\mu\nu\lambda\rho} e^{\mu}_1 e^{\nu}_3  e^{\lambda}_1 e^{\rho}_3 \over
\vert e_1 \wedge e_3\vert^2}= -  \ln^2 \lambda_2 < 0
\\
K_{23} = {R_{\mu\nu\lambda\rho} e^{\mu}_2 e^{\nu}_3  e^{\lambda}_2 e^{\rho}_3 \over
\vert e_2 \wedge e_3 \vert^2}= -  \ln^2 \lambda_1 <0. \nn
\eea
It follows for the above equations that the geodesic
flow is exponentially unstable on the planes (1,3) and (2,3)
and is stable in the plane (1,2). This behaviour is dual to the flow (\ref{velo}) which
is unstable in (1,2) plane and is stable in (1,3) and (2,3) planes.
 The scalar curvature is
\be
R= R_{\mu\nu\lambda\rho} g^{\mu\lambda}g ^{\nu\rho} = 2(K_{12} +K_{13}+K_{23}) =   - 2 
\ln^2 \lambda_1 = -2 h(T)^2 ,
\ee
where $h(T)$ is the entropy of the automorphism T (\ref{arno}).

\vfill


\begin{thebibliography}{99}
 
 

\bibitem{anosov}  D. V. Anosov, \emph{Geodesic flows on closed Riemannian manifolds with negative curvature},  Trudy Mat. Inst. Steklov., Vol. {\bf 90} (1967) 3 - 210

\bibitem{yer1986a}  G. Savvidy and N. Ter-Arutyunyan-Savvidy, \emph{ On the Monte Carlo simulation of physical systems}, J.Comput.Phys. {\bf 97} (1991) 566; Preprint EFI-865-16-86-YEREVAN, Jan. 1986. 13pp.

\bibitem{konstantin} K.Savvidy, \emph{The MIXMAX random number generator},
Comput.Phys.Commun. {\bf 196} (2015) 161; \url{http://dx.doi.org/10.1016/j.cpc.2015.06.003}



\bibitem{kolmo} A.N. Kolmogorov,  \emph{New metrical invariant of transitive dynamical 
systems and automorphisms of Lebesgue spaces}, 
Dokl. Acad. Nauk SSSR,  {\bf{119}} (1958) 861-865


\bibitem{kolmo1} A.N. Kolmogorov,  \emph{On the entropy per unit time as a metrical invariant
of automorphism}, 
Dokl. Acad. Nauk SSSR,  {\bf{124}} (1959) 754-755


\bibitem{Shannon}C.E.Shannon, \emph{The Mathematical Theory of Communications}, The Bell System Technical Journal, {\bf 27} (1948) 379-423 and 623-656

\bibitem{sinai3} Ya.G. Sinai,  \emph{On the Notion of Entropy of a Dynamical System}, Doklady of Russian Academy of Sciences, {\bf{124}}  (1959)  768-771.

\bibitem{Savvidy:2015jva}
K.~Savvidy and G.~Savvidy,
\emph{Spectrum and Entropy of C-systems. MIXMAX random number generator,}
Chaos Solitons Fractals \textbf{91} (2016), 33-38
doi:10.1016/j.chaos.2016.05.003
[arXiv:1510.06274 [math.DS]].
 
 
\bibitem{Savvidy:2015ida}
G.~K.~Savvidy,
 \emph{Anosov C-systems and random number generators,}
Theor. Math. Phys. \textbf{188} (2016) no.2, 1155-1171
doi:10.1134/S004057791608002X
[arXiv:1507.06348 [hep-th]].

\bibitem{Koopman} Koopman,\emph{Hamiltonian Systems and Transformations in Hilbert Space},
Proc. Nat. Acad. Sci.  {\bf 17} (1931) 315


 \bibitem{Rokhlin0} V. A. Rokhlin, \emph{Selected topics from the metric theory of dynamical systems}, Uspekhi Mat. Nauk, {\bf 4} (1949), 57-128    


 \bibitem{rokhlin1} V. A. Rokhlin, \emph{Metric properties of endomorphisms of compact commutative groups}, Izv. Akad. Nauk SSSR Ser. Mat., {\bf{28}} (1964) 867- 874   


\bibitem{rokhlin} V.A. Rokhlin, \emph{On the endomorphisms of compact commutative groups}, 
Izv. Akad. Nauk, {\bf{13}} (1949) 329 

\bibitem{rokhlin2}V.A. Rokhlin, \emph{On the entropy of automorphisms of compact 
commutative groups},  Teor.Ver. i Pril.,  {\bf3} (1961)  351


\bibitem{Savvidy:2018ygo}
G.~Savvidy and K.~Savvidy,
 \emph{Exponential decay of correlations functions in MIXMAX generator of pseudorandom numbers,}
Chaos Solitons Fractals {\bf 107} (2018), 244-250
doi:10.1016/j.chaos.2018.01.007 

\bibitem{Martirosyan:2018bjq}
N.~Martirosyan, K.~Savvidy and G.~Savvidy, \emph{Spectral Test of the MIXMAX Random Number Generators,} Chaos Solitons Fractals Solitons and Fractals: the interdisciplinary journal of Nonlinear Science {\bf 118} (2019), 242-248
doi:10.1016/j.chaos.2018.11.024
[arXiv:1806.05243 [nlin.CD]].



\bibitem{sinai2}Ya. G. Sinai, \emph{Markov partitions and C-diffeomorphisms}, Funkcional. Anal, i Prilozen. {\bf 2} (1968),64-89; Functional Anal. Appl. {\bf 2} (1968) 61-82.  

\bibitem{sinai4}Ya. G. Sinai, \emph{Proceedings of the International Congress 
of Mathematicians}, Uppsala (1963) 540-559.  

 
\bibitem{gines}A.~L.~Gines, \emph{Metrical properties of the endomorphisms on m-dimensional 
torus},  Dokl. Acad. Nauk SSSR,  {\bf 138} (1961) 991-993

\bibitem{Tsallis}  C. Tsallis, \emph{Possible Generalisation of Boltzmann-Gibbs Statistics},
 Journal of Statistical Physics, {\bf 52}  (1988) 479-488
 
\bibitem{Tsallis1} C. Tsallis, M. Gell-Mann and Y. Sato, \emph{Extensivity and entropy production},
Europhysics news {\bf 36} (2005) 186 -189?  
     
\bibitem{hepforge} K.Savvidy, MIXMAX code C/C++\\
HEPFORGE.ORG, {http://mixmax.hepforge.org}, \\
\url{http://www.inp.demokritos.gr/~savvidy/mixmax.php}
 
\bibitem{CLHEP} The foundation library CLHEP\\
\url{http://proj-clhep.web.cern.ch/proj-clhep/}\\
\url{https://gitlab.cern.ch/CLHEP/CLHEP/-/blob/develop/Random/Random/MixMaxRng.h}

\bibitem{Geant4} Geant4. Concurrent and Distributed MC toolkit, \\
\url{http://geant4.web.cern.ch}

\bibitem{root}  ROOT. Data analysis framework, \\
\url{https://root.cern.ch/doc/master/classROOT_1_1Math_1_1MixMaxEngine.html}\\
\url{https://root.cern.ch/doc/master/classTRandom.html}\\
 \url{https://root.cern.ch/doc/master/mixmax_8h_source.html}

\bibitem{CMSrunII}V. Ivanchenko and S. Banerjee,\emph{Upgrade of CMS Full Simulation for Run 2},
EPJ Web of Conferences {\bf 214} (2019)  02012,
\url{https://doi.org/10.1051/epjconf/201921402012}


\bibitem{CMS}V. Ivanchenko,
\url{https://indico.cern.ch/event/731433/contributions/3015654/attachments/1680131/2698971/CMSsim.pdf}\\
\url{https://indico.cern.ch/event/587955/contributions/2937635/attachments/1679273/2706817/PosterCMS_SIM_v4.pdf}\\

\bibitem{PYTHIA} T. Sj\"ostrand et al, \emph{An Introduction to PYTHIA 8.2} 
Comput. Phys.Commun. 191 (2015) 159 [arXiv:1410.3012 [hep-ph]], 
\url{http://home.thep.lu.se/~torbjorn/pythia83html/Welcome.html}

\bibitem{GSL} GSL-GNU Scientific Library, \emph{GNU Operating System}, \url{https://www.gnu.org/software/gsl/}


\bibitem{hedlund}G.Hedlund, \emph{The dynamics of geodesic flow}, 
Bull.Am.Math.Soc. {\bf 45} (1939) 241-246

\bibitem{hopf} E.Hopf. \emph{ Statistik der L\"osungen   
geod\"atischer Probleme vom unstabilen 
Typus.  II.}  Math.Ann. {\bf 117} (1940) 590-608



\bibitem{RCARRY}
G. Marsaglia and A. Zaman,   Ann.Appl.Probab. {\bf 1} (1991) 462-480.

\bibitem{MT}
M. Matsumoto and T. Nishimura, \emph{Mersenne Twister: A 623-dimensionally equidistributed uniform pseudorandom number generator}, ACM Trans. on Modelling and Computer Simulation {\bf 8} (1998) 3-30;  DOI:10.1145/272991.272995




 
\end{thebibliography}
\end{document}